\documentclass[11pt,a4paper]{article}
\usepackage[utf8]{inputenc}
\usepackage[T2A]{fontenc}
\usepackage[margin=2cm]{geometry}
\usepackage{graphicx}

\usepackage{amsmath,amsthm,amssymb}

\newtheorem{teo}{Theorem}

\newtheorem{lem}{Lemma}

\theoremstyle{definition}

\newtheorem{rem}{Remark}

\begin{document}

\title{On the Picard-Lindel\"{o}f Argument and the Banach-Caccioppoli Contraction Mapping Principle}

\def\firstname#1{#1}
\def\surname#1{#1}

\author{\firstname{Alexander\,I.}~\surname{Bufetov},
	\firstname{Ilya\,I.}~\surname{Zavolokin}}

\date{}
\maketitle

\begin{abstract}
	The aim of this note is to present the simple observation that a slight refinement of the Contraction Mapping Principle allows one to recover the precise convergence rate in the Picard-Lindel{\"o}f Theorem.\\
	\textbf{keywords:} Picard-Lindel{\"o}f Argument, Banach-Caccioppoli Contraction Mapping Principle, existence and uniqueness of solutions, ordinary differential equation
	
\end{abstract}
\maketitle


\begin{flushright}
	\textit{To our teacher Valery Vasilievich Kozlov, with admiration and love}				
\end{flushright}

\section{Introduction}
Ivan Georgievich Petrovsky in his 1939 textbook \cite{Petrovskiy1} gives the following three proofs of the existence and uniqueness of solutions for ordinary differential equations:
\begin{enumerate}
	\item Existence is proved by the Euler method of piecewise-linear approximation of solutions. Existence of two distinct solutions with the same initial condition is then shown to contradict the Lipschitz condition;
	\item the unique solution is found by iteration of an integral operator (the Picard-Lindel\"{o}f argument);
	\item existence and uniqueness are proved using the Banach-Caccioppoli Contraction Mapping Principle, to which I. G. Petrovsky himself gives the name ``Caccioppoli-Tikhonov Theorem$"$\footnote{In the 1949 edition, the Contraction Mapping Principle becomes the ``Tikhonov-Caccioppoli Theorem$"$,\, while in the posthumous edition \cite{Petrovskiy2}, the ``Banach Theorem$"$.}.
\end{enumerate}
\quad In his first proof Petrovsky uses a compactness argument based on the Arzel\`a-Ascoli Theorem, which, therefore, only yields the existence of solutions and does not give any approximation to solutions; note, however, that the Euler method, used in the proof of the Peano Existence Theorem, can be shown to yield a sequence of piecewise linear approximations converging at a polynomial rate. The Contraction Mapping Principle yields an exponential rate. The Picard-Lindel\"{o}f argument yields an estimate of the form $O(\alpha^n/n!)$ for the uniform norm of the difference between the $n$-th Picard approximation and the solution. An application of the Contraction Mapping Principle thus gives a slower convergence rate than the original method of Picard. Indeed Philip Hartman in his classical textbook \cite{Hartman} gives the Picard-Lindel\"{o}f argument only, not mentioning contraction maps at all, and so does Vyacheslav Vasilyevich Stepanov in  \cite{Stepanov}. Lev Semenovich Pontryagin in his ``Ordinary Differential Equations$"$ \cite{Pontryagin} gives both the Picard-Lindel\"{o}f argument and the proof based on the Contraction Mapping Principle, with a geometric rate.\\
\quad The aim of this note is to present the simple observation that a slight refinement of the Contraction Mapping Principle gives the rate of $O(\alpha^n/n!)$; namely, we consider a decreasing chain of subspaces, and a mapping that takes each of the subspaces of the chain into the next one and contracts more and more on each subsequent subspace.\\
\quad For instance, let $H$ be a complete metric space with metric $d$, let $H_n\subset H$ be a decreasing chain of closed subspaces:
$$
H_0 = H \supset H_1 \supset H_2 \ldots \supset H_n \supset \ldots,
$$
let $\varkappa_n\in (0, 1)$ and let $T: H \to H$ be a map satisfying the assumptions:
\begin{enumerate}
	\item $TH_n\subset H_{n+1},$ for all $n\geqslant 0$;
	\item $d(Tx, Ty) \leqslant \varkappa_n d(x, y)$, for all $x, y\in H_n,\, n\geqslant 0$.
\end{enumerate}
Then the map $T$ has a unique fixed point $p$, and for all $x\in H$ we have the inequality:
\[d(T^n x, p) \leqslant \varkappa_1\ldots \varkappa_{n+1}d(x, p),\, n\geqslant 1.\]
A slight modification of this simple scheme -- a different metric must be chosen on each $H_n$ -- provides a proof of the Existence and Uniqueness Theorem for  ordinary differential equations and yields the precise rate of convergence. We conclude this note with a few brief remarks of the history of proofs of the Existence and Uniqueness Theorems from Euler and Cauchy to Picard, Lindel\"{o}f, Banach and Caccioppoli.

\section{A contraction mapping on a chain of metric spaces}
Let $\alpha_n>0$ be positive constants and let $(H_j, d_j)_{j=0}^\infty$ be a decreasing chain of complete metric spaces:
$$H_0 \supset H_1 \supset H_2 \supset \ldots$$
such that we have
\begin{equation}\tag{1}
	d_j (x, y) \leqslant \alpha_{j+1} d_{j+1} (x, y),\; \text{for all}\, x, y\in H_{j+1},\label{metrics}
\end{equation}
Let $\varkappa_j>0$ be constants such that
\begin{equation}
	\limsup\limits_j \alpha_j \varkappa_j < 1\tag{2}\label{kappaineq}
\end{equation}
and let $P:H_0 \to H_0$ be a map with following properties:
\begin{equation}
	P H_j \subset H_{j+1},\; \text{for all}\, j=0,1,\ldots \tag{3}\label{Scale spaces}
\end{equation}
\begin{equation}
	d_{j+1}(P(x), P(y)) \leqslant \varkappa_{j+1} d_j(x, y),\, \forall x, y\in H_j,\; j=0,1,\ldots \tag{4}\label{Contractions}
\end{equation}
\begin{rem}
	In the application to the Picard-Lindel{\"o}f argument, the constant $\alpha_1 = \alpha_2 = \ldots = \alpha$ gives the length of the interval on which the solution is defined.
\end{rem}
Set \begin{displaymath} 
	S = \bigcap_j H_j.
\end{displaymath}
\begin{lem}\label{lem:1}
	The map $P$ has a unique fixed point $x_\infty\in S$. Moreover, for all $x\in H$ and for all $n\geqslant j = 0, 1, \ldots$ the sequence $(x_m)_m$ defined by the formula:
	\[x_0 = x,\; x_m = P^{\circ m} x,\, m\geqslant 1\]
	satisfies the relations
	\[d_j(x_{\infty}, x_n)\leqslant C \alpha_n\varkappa_{n}\alpha_{n-1}\varkappa_{n-1} \ldots \alpha_{j+1}\varkappa_{j+1}d_j(x_{1+j}, x_j),\]
	where
	\[C = \sup_n \sum_{m=0}^\infty \alpha_{n+m}\varkappa_{n+m} \alpha_{n+m-1}\varkappa_{n+m-1} \ldots \alpha_{n+1}\varkappa_{n+1}\tag{5}\label{constant series}.\]
\end{lem}
\begin{rem}
	The series (\ref{constant series}) converges by (\ref{kappaineq}).
\end{rem}
\begin{rem}
	In particular, if $\varkappa_n = L/n$, and $\alpha_1 = \alpha_2 = \ldots = \alpha$ then setting $j=0$, we obtain:
	\[d_0(x_{\infty}, x_n)\leqslant e^{\alpha L} \frac{(L \alpha)^{n}}{n!}d_0(x_{1}, x_0).\]
\end{rem}
\begin{rem}
	If $d_j = d_0\Bigl|_{H_j},\, j=1,\ldots,$ then $\alpha = 1$ and the proof is immediate --- but, as we will see later, it is insufficient for the Picard-Lindel\"{o}f theorem, cf. Remark~\ref{rem:5} below.
\end{rem}
\begin{proof}
	We will show that there exists a point $x_\infty$ from $S$ the intersection of the spaces $H_j,\, j\geqslant 0$, which is attracting in all metrics $d_j,\,j\geqslant0:$
	$$\exists x_\infty \in S: \forall j=0,1,\ldots\;\forall x\in H_j: P^{\circ n} x\xrightarrow{d_j} x_\infty.$$
	Take $x\in H$ and denote $x_m:= P^{\circ m} x.$ Then for $n\geqslant j =0, 1, \ldots$ we have the estimates
	\[d_j(x_{n+1}, x_n)\leqslant \alpha_n\alpha_{n-1}\ldots\alpha_{j+1} d_n(x_{n+1}, x_n) \leqslant \alpha_n\alpha_{n-1}\ldots\alpha_{j+1}\varkappa_n\varkappa_{n-1} \ldots \varkappa_{j+1} d_j(x_{1+j}, x_j),\]
	where the first inequality was obtained using (\ref{metrics}) and the second obtained by the properties \ref{Scale spaces}, \ref{Contractions}) of the operator $P$. For $m\geqslant 0$, we obtain the estimate
	\begin{gather*}
		d_j(x_{n+m}, x_n)\leqslant d_j(x_{n+m}, x_{n+m-1}) + \ldots + d_j(x_{n+1}, x_{n}) \leqslant \tag{6}\label{rate}\\ 
		(\alpha_{n+m}\alpha_{n+m-1}\ldots\alpha_{j+1}\varkappa_{n+m}\varkappa_{n+m-1} \ldots \varkappa_{j+1} + \ldots + \alpha_{n}\alpha_{n-1}\ldots\alpha_{j+1}\varkappa_{n}\varkappa_{n-1} \ldots \varkappa_{j+1}) d_j(x_{1+j}, x_j) = \\
		\alpha_n\varkappa_{n} \alpha_{n-1}\varkappa_{n-1} \ldots \alpha_{j+1}\varkappa_{j+1} (\alpha_{n+m}\varkappa_{n+m}\alpha_{n+m-1}\varkappa_{n+m-1} \ldots \alpha_{n+1}\varkappa_{n+1} + \ldots + 1) d_j(x_{1+j}, x_j)
	\end{gather*}
	Since $$\limsup \alpha_j\varkappa_j < 1,\,$$ it follows that the series $\sum_{m=0}^\infty \alpha_{n+m}\varkappa_{n+m} \alpha_{n+m-1}\varkappa_{n+m-1} \ldots \alpha_{n+1}\varkappa_{n+1}$ converges for all $n$, and we set
	\[S(n) = \sum_{m=0}^\infty \alpha_{n+m}\varkappa_{n+m} \alpha_{n+m-1}\varkappa_{n+m-1} \ldots \alpha_{n+1}\varkappa_{n+1}.\]
	By assumption \ref{kappaineq}, for large enough $n$ we have:
	\[\alpha_{n+m}\varkappa_{n+m} \alpha_{n+m-1}\varkappa_{n+m-1} \ldots \alpha_{n+1}\varkappa_{n+1} + \ldots + 1 < (\limsup \alpha_j\varkappa_j)^m + \ldots + 1 < \infty;\]
	whence 
	$$\sup_n S(n) = C<\infty.$$
	Therefore for any $j\in \mathbb N$, $(x_n)_n$ is a Cauchy sequence in the metric space $(H_j, d_j)$, and we have 
	$$x_n \xrightarrow{d_j} x_\infty^j \in H_j.$$
	
	We now check the uniqueness of the fixed point. If $z \in H_0$ is such that $P(z) = z$ then for all $j\geqslant0$ we have:
	\[d_0(x_\infty, z) = d_0(P x_\infty, P z) < \alpha_j \varkappa_j\alpha_{j-1}\varkappa_{j-1}\ldots \alpha_1\varkappa_1 d_0(x_\infty, z),\]
	and $d_0(x_\infty, z) =0.$\\
	We now estimate the rate of convergence. Passing to the limits as $m\to\infty$ in (\ref{rate}), we obtain the inequality:
	\[d_j(x_{\infty}, x_n)\leqslant C \alpha_{n}\varkappa_{n} \alpha_{n-1}\varkappa_{n-1} \ldots \alpha_{j+1}\varkappa_{j+1}d_j(x_{1+j}, x_j).\]
	The desired rate of convergence is established. Lemma \ref{lem:1} is proved completely.
\end{proof}
\section{A proof of the Picard-Lindel\"of Theorem}
\subsection{Preliminary remarks}
We now apply the Lemma to the proof of the Picard-Lindel\"{o}f Theorem. 

Endow the space $\mathbb R^d,\,d\in\mathbb N$ with the norm 
$$||v|| = \sqrt{\sum_{i=1}^d v_i^2}.$$
Let $t_0\in \mathbb R,\, y_0\in\mathbb R^d,\, a, b>0.$ Consider the space-time parallelepiped:
\[R=\{(t, y): |t-t_0|\leqslant a,\, ||y-y_0||\leqslant b\}\]
and let $f: R \to \mathbb R^d,$ be a continuous function satisfies the Lipshitz condition in y:
\[|f(t, y_1) - f(t, y_2)|\leqslant L||y_1 - y_2||,\]
and $|f(x)|\leqslant M$ on $R$. Let $\alpha = \min(a, b/M)$.
\begin{rem}\label{rem:5}
	As in the classical Picard-Lindel{\"o}f argument, the constant $\alpha$ gives the length of the interval on which the solution is defined.
\end{rem}

Introduce the complete metric space $H_0$ of all continuous maps
$$y:[t_0-\alpha, t_0 + \alpha]\to \mathbb R^d$$
satisfying the assumptions $||y(t)-y_0||\leqslant b$, and $y(t_0)=y_0$ with induced uniform metric 
$$d_0(x, y) = \max_{t\in[t_0-\alpha, t_0 + \alpha]} ||x(t)-y(t)||.$$

Introduce the Picard operator $P:H_0 \to H_0$ by the formula
\begin{equation}
	(Py)(t) = y_0 + \int_{t_0}^t f(s, y(s))ds \tag{7}\label{Picard operator}
\end{equation}
\begin{rem}
	Suppose $y$ is a fixed point for Picard operator. By definition we have
	\[y(t_0) = y_0,\]
	while differentiating (\ref{Picard operator}) yields
	\begin{equation}
		\frac{d}{dt}y(t) = f(t, y(t)),\tag{8}\label{Cauchy problem}
	\end{equation}
	whence $y$ is the solution of Cauchy problem for the equation (\ref{Cauchy problem}).
\end{rem}

\begin{teo}[Picard-Lindel\"{o}f]
	There exists a unique continuously differentiable function $y^\infty \in H_0$
	such that $Py^\infty = y^\infty$ and for any continuous function
	$$y: \{|t-t_0|\leqslant\alpha\} \to \{||y-y_0||\leqslant b\}$$
	with $y(t_0) = y_0$ the sequence defined by the formula $y^n = P^{\circ n}y$ converges to $y^\infty\in H_0$, and satisfies the estimate
	\[||y^n - y^\infty||_{C[t_0-\alpha, t_0+\alpha]}\leqslant e^{\alpha L} \frac{(\alpha L)^n}{n!}M.\]
\end{teo}
\begin{rem}
	As mentioned above, we cannot simplify our presentation by using the induced metrics $d_j = d_0\Bigl|_{H_j}$. Indeed, for the Picard operator $P$ we have
	\begin{multline*}
		d(Px, Py) = \max_{t\in[t_0-\alpha,t_0+\alpha]} \left|\left|\int_{t_0}^t f(s, x(s)) - f(s, y(s)) ds\right|\right| \leqslant\\
		\leqslant\max_{t\in[t_0-\alpha,t_0+\alpha]} \int_{t_0}^t L||x(s) - y(s)|| |ds|
		\leqslant \int_{t_0}^{t_0+\alpha} L d(x, y)ds = \alpha L d(x, y).
	\end{multline*}
	Thus, we only obtain convergence in the interval $[t_0-\alpha, t_0+\alpha]$ with $\alpha L<1$ and only with geometric convergence rate.
\end{rem}
\begin{rem}
	The norm $||\cdot||$ in $\mathbb R^d$ is chosen for concreteness only; any other norm can be used instead. The argument works in general Banach spaces.
\end{rem}
	\subsection{A chain of subspaces}
	For $j\in\mathbb N$, let $H_j$ be the subspace of functions $y\in H_0$ satisfying the condition:
	\[\sup_{\substack{t\in[t_0-\alpha, t_0+\alpha]\\ t\ne t_0}} \frac{||Py(t) - y(t)||}{|t-t_0|^j} <\infty \tag{9}\label{sup-j}\]
	Let $C_j(f, y)$ stand for supremum in (\ref{sup-j}). Define the metric $d_j$ on $H_j$ by the formula:
	\[d_j(x, y) = \sup_{\substack{t\in[t_0-\alpha, t_0+\alpha]\\t\ne t_0}} \left|\left|\frac{x(t)-y(t)}{(t-t_0)^j}\right|\right|.\tag{10}\label{sup-j-1}\]
	We now check that the supremum in \eqref{sup-j-1} is finite. Define the following linear bounded operator $\mathcal K$ on the Banach space $C[t_0-\alpha, t_0+\alpha]:$
	\[\left(\mathcal{K}g\right)(t) = \int_{t_0}^t g(s)|ds|,\]
	for $g\in C[t_0-\alpha, t_0+\alpha]$. It is easy to verify, that for $n\in\mathbb N$ we have:
	\[\left(\mathcal{K}^n g\right)(t) = \int_{t_0}^t g(s)\frac{|t-s|^{n-1}}{(n-1)!}|ds|,\]
	and
	\begin{align*}
		(Id - L\mathcal{K})^{-1}g(t) = g(t) + \sum_{n=1}^\infty (L\mathcal K)^n g(t) = g(t) + \sum_{n=1}^\infty L^n \int_{t_0}^t g(s)\frac{|t-s|^{n-1}}{(n-1)!}|ds| =\\
		= g(t) + L \sum_{n=0}^\infty L^n \int_{t_0}^t g(s)\frac{|t-s|^n}{n!}|ds| = g(t) + L \int_{t_0}^t g(s) e^{L|t-s|}|ds|
	\end{align*}
	Write
	\[\left[R_L g\right](t) = g(t) + L \int_{t_0}^t g(s) e^{L|t-s|}|ds|,\]
	for all $g\in C[t_0-\alpha, t_0+\alpha].$ By definition, the operator $R_L$ is monotone in following sense:
	\[g(t)\leqslant h(t),\,\forall t\in[t_0-\alpha, t_0+\alpha] \implies R_L g(t) \leqslant R_L h(t),\,\forall t\in[t_0-\alpha, t_0+\alpha].\]
	For all $x, y \in H_j$, we have
	\[||x(t)-y(t)|| \leqslant ||x(t) - Px(t)|| + ||y(t) - Py(t)|| + ||Px(t)-Py(t)||.\]
	Since,
	\[||Px(t) - Py(t)|| \leqslant \int_{t_0}^t ||f(s, x(s)) - f(s, y(s))|| |ds| \leqslant L \int_{t_0}^t ||x(s) - y(s)|| |ds|,\]
	we obtain
	\[||x(t) - y(t)|| - L \int_{t_0}^t ||x(s) - y(s)|| |ds| = (Id - L\mathcal{K})||x-y||(t)\leqslant (C_j(f, x) + C_j(f, y))|t-t_0|^j.\]
	Therefore,
	\begin{multline*}
		||x-y||(t) = (Id - L\mathcal K)^{-1}(Id - L\mathcal K)||x-y||(t) \leqslant (C_j(f, x) + C_j(f, y))(Id - L\mathcal K)^{-1}|t-t_0|^j = \\
		=(C_j(f, x) + C_j(f, y))\left(|t-t_0|^j + \sum_{k=j+1}^\infty \frac{L^k |t-t_0|^k}{(j+1)(j+2)\ldots k}\right) =\\
		=(C_j(f, x) + C_j(f, y))\left(1 + L^j \sum_{k=1}^\infty \frac{L^k |t-t_0|^k}{(j+1)(j+2)\ldots (j+k)}\right)|t-t_0|^j\leqslant\\
		\leqslant (C_j(f, x) + C_j(f, y))\left(1 + L^j \sum_{k=1}^\infty \frac{L^k |t-t_0|^k}{k!}\right)|t-t_0|^j =\\
		=(C_j(f, x) + C_j(f, y))\left(1 + L^j (e^{L|t-t_0|}-1)\right)|t-t_0|^j \leqslant (C_j(f, x) + C_j(f, y))\left(1 + L^j (e^{L\alpha}-1)\right)|t-t_0|^j,
	\end{multline*}
	we finally arrive at the desired estimate
	\[d_j(x, y) < (C_j(f, x) + C_j(f, y))\left(1 + L^j (e^{L\alpha}-1)\right).\]
	By definition the metrics $d_j$ satisfy the inequality
	\[
		d_j(x, y) = \sup_{\substack{t\in[t_0-\alpha, t_0+\alpha]\\t\ne t_0}} \left|\left|\frac{x(t)-y(t)}{(t-t_0)^j}\right|\right| \leqslant\sup_{\substack{t\in[t_0-\alpha, t_0+\alpha]\\t\ne t_0}}\left( \left|\left|\frac{x(t)-y(t)}{(t-t_0)^{j+1}}\right|\right| |t-t_0|\right)\leqslant \alpha d_{j+1}(x, y),
	\]
	for all $x, y\in H_{j+1}$. 
	\subsection{Completeness of the spaces $H_j$}
	\begin{rem}
		One way of checking completeness of the space $H_j$ is to take $y_0, y_1, y_2 \in H_j,$ to set
		\[u_{10}(t) = \dfrac{y_1(t) - y_0(t)}{(t-t_0)^j},\quad u_{20}(t) = \dfrac{y_2(t) - y_0(t)}{(t-t_0)^j},\]
		and note that, be definition, we have $u_{10}, u_{20}\in C[t_0-\alpha, t_0+\alpha]$ and that
		\[d_j(y_1, y_2) = \max\limits_{t:|t-t_0|<\alpha} |u_{10}(t)-u_{20}(t)|.\]
		Completeness of Tchebycheff metric now implies completeness of the metric $d_j$. One can also proceed as follows.
	\end{rem}
	We check that metric spaces $(H_j, d_j)$ are complete. Indeed, consider a Cauchy sequence $(x_n)$. The inequality
	\[d_0(x_n, x_m) \leqslant \alpha^j d_j(x_n, x_m) \to 0\]
	implies that there exists a function $x\in C[t_0-\alpha, t_0+\alpha]$ such that
	\[d_0 (x_n, x) \to 0.\]
	From the inequality
	\[d_j(x_n, x_m) \geqslant \sup_{\substack{t\in[t_0-\alpha,t_0+\alpha]\\t\ne t_0}} \left(\frac{||x_n(t) - Px_n(t)||}{|t-t_0|^j} - \frac{||x_m(t) - Px_m(t)||}{|t-t_0|^j} -\frac{||Px_n(t) - Px_m(t)||}{|t-t_0|^j}\right),\]
	we obtain
	\[d_j(x_n, x_m) \geqslant C_j(f, x_n) - C_j(f, x_m) - d_j(Px_n, Px_m) \geqslant C_j(f, x_n) - C_j(f, x_m) - \alpha L d_j(x_n, x_m),\]
	whence
	\[C_j(f, x_n) - C_j(f, x_m) \leqslant (1+\alpha L)d_j(x_n, x_m).\]
	Similarly,
	\[C_j(f, x_m) - C_j(f, x_n) \leqslant (1+\alpha L)d_j(x_n, x_m).\]
	The sequence of real numbers $(C_j(f, x_n))_n$ is Cauchy sequence, and we let $C$ be its limit:
	\[\lim\limits_{n\to\infty} C_j(f, x_n) = C<\infty.\]
	Taking $n\to\infty$ in the inequality
	\[||Px_n(t) - x_n(t)|| \leqslant C_j(f, x_n)|t-t_0|^j,\]
	we obtain
	\[||Px(t) - x(t)|| \leqslant C |t-t_0|^j.\]
	Therefore, $x\in H_j$ and the metric space $H_j$ is complete.
	\subsection{Nested images}
	The Picard operator $P$ map $H_j$ into $H_{j+1}$, because for all $y\in H_j,\,j =0, 1, 2,\ldots$ we have
	\begin{multline*}
		||P(Py) - Py|| = \left|\left|\int_{t_0}^t (f(s, Py(s)) - f(s, y(s)))ds\right|\right| \leqslant\\
		\leqslant \int_{t_0}^t L ||Py(s) - y(s)|| |ds| \leqslant \frac{C_j(f, y) L}{j+1} |t-t_0|^{j+1}
	\end{multline*}
	and
	\[||Py - y_0|| = \left|\left|\int_{t_0}^t f(s, y(s)) ds\right|\right| \leqslant \int_{t_0}^t ||f(s, y(s))|| |ds|\leqslant \int_{t_0}^\alpha M ds= b\]
	\subsection{The rate of convergence}
	In order to apply Lemma~\ref{lem:1} it remains to prove that assumptions \eqref{kappaineq}, and \eqref{Contractions} hold. We have
	\begin{multline*}
		d_{j+1}(Px, Py) = \sup_{\substack{t\in[t_0-\alpha, t_0+\alpha]\\t\ne t_0}} \frac{1}{|t-t_0|^{j+1}}\left|\left|\int_{t_0}^t f(s, x(s)) - f(s, y(s))ds\right|\right| \leqslant\\
		\leqslant \sup_{\substack{t\in[t_0-\alpha, t_0+\alpha]\\t\ne t_0}} \frac{1}{|t-t_0|^{j+1}}\int_{t_0}^{t} L ||x(s) - y(s)|| |ds|\leqslant \\
		\leqslant \sup_{\substack{t\in[t_0-\alpha, t_0+\alpha]\\t\ne t_0}} \frac{1}{|t-t_0|^{j+1}}\int_{t_0}^{t} L \sup_{s\ne t_0}\left(\frac{||x(s)-y(s)||}{|s - t_0|^j}\right) |s-t_0|^j |ds| = \frac{L}{j+1} d_{j}(x, y).
	\end{multline*}
	In particular, we have
	$$\varkappa_j = \frac{L}{j}\to 0$$
	and the assumptions of Lemma~\ref{lem:1} are verified. Set $y^n = P^{\circ n} y$. Using Lemma~\ref{lem:1}, we obtain $y^n \xrightarrow{d_0} y^\infty.$ Since
	$$y^\infty\in \bigcap_j H_j,$$
	it follows that
	\[y^\infty(t) = (Py^\infty)(t) = y_0 + \int_{t_0}^t f(s, y^\infty(s))ds,\]
	and we obtain the desired rate of convergence
	\[d_0(y^n, y^\infty) \leqslant e^{\alpha L} \frac{(\alpha L)^n}{n!}M.\] 
\begin{rem}
	One could define the subspaces $H_j$ simply as images of the space $H_0$ under the Picard operator. We prefer an intrinsic characterization: the space $H_j$ consists of those functions for which the Picard operator does not change its $j$-th jet.
\end{rem}
\section{Differential equations with complex time}
\subsection{Preliminary remarks}
The same argument works in complex time. Endow the space $\mathbb C^d,\,d\in\mathbb N$, with the norm 
$$||v|| = \max_{j=1,\ldots, d} |v_j|.$$
Let $t_0\in \mathbb C,\, z_0\in\mathbb C^d,\, a, b>0.$ Consider the complex space-time parallelepiped:
\[R=\{(t, y): |t-t_0|\leqslant a,\, ||z-z_0||\leqslant b\}\]
and let $F: R \to \mathbb C^d,$ be a holomorphic function satisfying the Lipshitz condition in $z$:
\[||F(t, z_1) - F(t, z_2)||\leqslant L||z_1 - z_2||,\,\text{for all}\, (t, z_1), (t, z_2)\in R\]
and $|F(t, z)|\leqslant M$ on $R$. Let $\alpha = \min(a, b/M)$, let $U = \{|t-t_0|<\alpha\}$ be a disc, and $V=\{||z-z_0||< b\}$ be a polydisc. 
Introduce the metric space $H_0$ of all, holomorphic in $U$, continuous in $\overline{U}$
$$z:\overline{U}\to \overline{V}$$
maps and
satisfying the assumptions $||z(t)-z_0||\leqslant b$, and $z(t_0)=z_0$. The space $H_0$ is endowed with the uniform metric 
$$d_0(z, w) = \max_{t\in \overline{U}}||z(t)-w(t)||.$$
Introduce the Picard operator $P:H_0 \to H_0$ by the formula
\[(Pz)(t) = z_0 + \int_{t_0}^t F(s, z(s))ds,\]
where the integration takes place over an internal joining $t_0$ to $t$.
\begin{teo}[The Picard-Lindel\"of Theorem in complex time]\label{PL_complex}
	There exists a unique holomorphic function
	$$z^\infty: U \to V$$
	such that $Pz^\infty = z^\infty$. For any holomorphic function
	$$z: U \to V$$
	with $z(t_0) = z_0$ the sequence defined by the formula $z^n = P^{\circ}z$ uniformly converges to $z^\infty$, and satisfies the estimate
	\[||z^n - z^\infty||_{C[\overline{U}]}\leqslant e^{\alpha L} \frac{(\alpha L)^n}{n!}M.\]
\end{teo}

\begin{rem}
	A proof of the Existence and Uniqueness Theorem in complex time based on the Contraction Mapping Principle can be found, for instance, in \cite{Ilyashenko-Yakovenko}.
\end{rem}

	\subsection{A chain of subspaces}
	Consider the metric space $(H_0, d_0)$. Recalling  that, by the Weierstrass Theorem, the limit of a uniformly convergent sequence of holomorphic functions is itself holomorphic, we see that the space $(H_0, d_0)$ is complete.

	For $j\in\mathbb N$, let $H_j$ be the subspace of functions $z\in H_0$ satisfying the condition:
	\[\sup_{\substack{t\in\overline{U}\\ t\ne t_0}} \frac{||Pz(t) - z(t)||}{|t-t_0|^j} <\infty \tag{13}\label{sup-j complex}\]
	Let $C_j(F, z)$ stand for supremum in (\ref{sup-j complex}). Define the metric $d_j$ on $H_j$ by the formula:
	\[d_j(z, w) = \sup_{\substack{t\in\overline{U}\\t\ne t_0}} \left|\left|\frac{z(t)-w(t)}{(t-t_0)^j}\right|\right|.\tag{14}\label{sup-j complex-1}\]

	We check that the supremum in \eqref{sup-j complex-1} is finite for all $z, w\in H_j$. Define the linear bounded operator $\mathcal K$ on the Banach space of real continuous functions $g\in C(\overline{U})$ by the formula
	\[\left(\mathcal{K}g\right)(t) = \int_{t_0}^t g(s)|ds|.\]
	Again, the integration takes place over the internal joining $t$ and $t_0$. For $n\in\mathbb N$ we have
	\[\left(\mathcal{K}^n g\right)(t) = \int_{t_0}^t g(s)\frac{|t-s|^{n-1}}{(n-1)!}|ds|,\]
	and
	\begin{align*}
		(Id - L\mathcal{K})^{-1}g(t) = g(t) + \sum_{n=1}^\infty (L\mathcal K)^n g(t) = g(t) + \sum_{n=1}^\infty L^n \int_{t_0}^t g(s)\frac{|t-s|^{n-1}}{(n-1)!}|ds| =\\
		= g(t) + L \sum_{n=0}^\infty L^n \int_{t_0}^t g(s)\frac{|t-s|^n}{n!}|ds| = g(t) + L \int_{t_0}^t g(s) e^{L|t-s|}|ds|
	\end{align*}
	Write
	\[\left[R_L g\right](t) = g(t) + L \int_{t_0}^t g(s) e^{L|t-s|}|ds|.\]
	For $g\in C(\overline{U}).$ By definition the operator $R_L$ is monotone in the following sense:
	\[g(t)\leqslant h(t),\,\forall t\in\overline{U} \implies R_L g(t) \leqslant R_L h(t),\,\forall t\in\overline{U}.\]
	For all $z, w \in H_j$, we have
	\[||z(t)-w(t)|| \leqslant ||z(t) - Pz(t)|| + ||w(t) - Pw(t)|| + ||Pz(t)-Pw(t)||.\]
	By definition of the Picard operator, we have
	\[||Pz(t) - Pw(t)|| \leqslant \int_{t_0}^t ||F(s, z(s)) - F(s, w(s))|| |ds| \leqslant L \int_{t_0}^t ||z(s) - w(s)|| |ds|,\]
	whence
	\[||z(t) - w(t)|| - L \int_{t_0}^t ||z(s) - w(s)|| |ds| = (Id - L\mathcal{K})||z - w||(t)\leqslant (C_j(F, z) + C_j(F, w))|t-t_0|^j.\]
	Therefore,
	\begin{multline*}
		||z - w||(t) = (Id - L\mathcal K)^{-1}(Id - L\mathcal K)||z-w||(t) \leqslant (C_j(F, z) + C_j(F, w))(Id - L\mathcal K)^{-1}|t-t_0|^j =\\
		= (C_j(F, z) + C_j(F, w))\left(|t-t_0|^j + \sum_{k=j+1}^\infty \frac{L^k |t-t_0|^k}{(j+1)(j+2)\ldots k}\right)=\\
		=(C_j(F, z) + C_j(F, w))\left(1 + L^j \sum_{k=1}^\infty \frac{L^k |t-t_0|^k}{(j+1)(j+2)\ldots (j+k)}\right)|t-t_0|^j\leqslant\\
		\leqslant (C_j(F, z) + C_j(F, w))\left(1 + L^j \sum_{k=1}^\infty \frac{L^k |t-t_0|^k}{k!}\right)|t-t_0|^j =\\
		=(C_j(F, z) + C_j(F, w))\left(1 + L^j (e^{L|t-t_0|}-1)\right)|t-t_0|^j \leqslant (C_j(F, z) + C_j(F, w))\left(1 + L^j (e^{L\alpha}-1)\right)|t-t_0|^j,
	\end{multline*}
	and
	\[d_j(z, w) < (C_j(F, z) + C_j(F, w))(1+ L^j (e^{L\alpha}-1)).\]
	The metric $d_j$ is thus well-defined. By definition the metrics $d_j$ satisfy the inequality
	\[
		d_j(z, w) = \sup_{\substack{t\in\overline{U}\\t\ne t_0}} \left|\left|\frac{z(t)-w(t)}{(t-t_0)^j}\right|\right| \leqslant\sup_{\substack{t\in\overline{U}\\t\ne t_0}}\left( \left|\left|\frac{z(t)-w(t)}{(t-t_0)^{j+1}}\right|\right| |t-t_0|\right)\leqslant \alpha d_{j+1}(z, w),
	\]
	for all $z, w\in H_{j+1}$. 
	\subsection{Completeness of the spaces $H_j$}
	\begin{rem}
		One way of checking completeness of the space $H_j$ is to take $y_0, y_1, y_2 \in H_j,$ to set
		\[u_{10}(t) = \dfrac{y_1(t) - y_0(t)}{(t-t_0)^j},\quad u_{20}(t) = \dfrac{y_2(t) - y_0(t)}{(t-t_0)^j},\]
		and note that, be definition, we have $u_{10}, u_{20}\in C(\overline U)$ and that
		\[d_j(y_1, y_2) = \max\limits_{t\in\overline U} |u_{10}(t)-u_{20}(t)|.\]
		Completeness of Tchebycheff metric now implies completeness of the metric $d_j$. One can also proceed as follows.
	\end{rem}
	
	We check that the metric spaces $(H_j, d_j)$ are complete. Indeed, consider a Cauchy sequence $(z_n)_n$. The inequality
	\[d_0(z_n, z_m) \leqslant \alpha^j d_j(z_n, z_m) \to 0\]
	implies that there exists a function $z\in C[\overline{U}]$ such that $z$ is holomorphic in restriction to $U$ and such that
	\[d_0 (z_n, z) \to 0.\]
	From the inequality
	\[d_j(z_n, z_m) \geqslant \sup_{\substack{t\in\overline{U}\\t\ne t_0}} \left(\frac{||z_n(t) - Pz_n(t)||}{|t-t_0|^j} - \frac{||z_m(t) - Pz_m(t)||}{|t-t_0|^j} -\frac{||Pz_n(t) - Pz_m(t)||}{|t-t_0|^j}\right),\]
	we obtain
	\[d_j(z_n, z_m) \geqslant C_j(F, z_n) - C_j(F, z_m) - d_j(Pz_n, Pz_m) \geqslant C_j(F, z_n) - C_j(F, z_m) - \alpha L d_j(z_n, z_m),\]
	whence
	\[C_j(F, z_n) - C_j(F, z_m) \leqslant (1+\alpha L)d_j(z_n, z_m).\]
	Similarly,
	\[C_j(F, z_m) - C_j(F, z_n) \leqslant (1+\alpha L)d_j(z_n, z_m).\]
	Therefore, the sequence $(C_j(F, z_n))_n$ is a Cauchy sequence of real numbers. Set
	\[C = \lim\limits_{n\to\infty} C_j(F, z_n)<\infty.\]
	Taking $n\to\infty$ in the inequality
	\[||Pz_n(t) - z_n(t)|| \leqslant C_j(F, z_n)|t-t_0|^j,\]
	we obtain
	\[||Pz(t) - z(t)|| \leqslant C |t-t_0|^j.\]
	
	Therefore, $z\in H_j$ and the metric space $H_j$ is complete. 
	\subsection{Nested images}
	The Picard operator $P$ maps $H_j$ into $H_{j+1}$, because for all $z\in H_j,\,j =0, 1, 2,\ldots$ we have
	\begin{multline*}
		||P(Pz) - Pz|| = \left|\left|\int_{t_0}^t (F(s, Pz(s)) - F(s, z(s)))ds\right|\right| \leqslant\\
		\leqslant \int_{t_0}^t L ||Pz(s) - z(s)|| |ds|\leqslant \frac{C_j(f, z) L}{j+1} |t-t_0|^{j+1}
	\end{multline*}
	and
	\[||Pz - z_0|| = \left|\left|\int_{t_0}^t F(s, z(s)) ds\right|\right| \leqslant \int_{t_0}^t ||F(s, z(s))|| |ds|\leqslant \int_{t_0}^\alpha M ds= b\]
	\subsection{The rate of convergence}
	In order to apply Lemma~\ref{lem:1} it remains to prove that its assumptions \eqref{kappaineq}, and \eqref{Contractions} hold. We have
	\begin{multline*}
		d_{j+1}(Pz, Pw) = \sup_{\substack{t\in\overline{U}\\t\ne t_0}} \frac{1}{|t-t_0|^{j+1}}\left|\left|\int_{t_0}^t F(s, z(s)) - F(s, w(s))ds\right|\right| \leqslant\\
		\leqslant \sup_{t}\frac{1}{|t-t_0|^{j+1}}\int_{t_0}^{t} L ||z(s) - w(s)|| |ds|\leqslant \\
		\leqslant \frac{1}{\alpha^{j+1}}\int_{t_0}^{t_0+\alpha} L \sup_s\left(\frac{||z(s)-w(s)||}{|s - t_0|^j}\right) |s-t_0|^j |ds| = \frac{L}{j+1} d_{j}(z(s), w(s)).
	\end{multline*}
	In particular, we have
	$$\varkappa_j = \frac{L}{j}\to 0.$$
	Set $z^n = P^{\circ n} z$. Using Lemma~\ref{lem:1}, we obtain $z^n \xrightarrow{d_0} z^\infty.$ Since
	$$z^\infty\in \bigcap_j H_j,$$
	it follows that
	\[z^\infty(t) = (Pz^\infty)(t) = z_0 + \int_{t_0}^t F(s, z^\infty(s))ds,\]
	and we obtain the desired rate of convergence
	\[d_0(z^n, z^\infty) \leqslant e^{\alpha L} \frac{(\alpha L)^n}{n!}M.\] 
	Theorem \ref{PL_complex} is proved completely.

	
\section{Historical remarks}
Gottfried Wilhelm Leibniz \cite{Leibniz} proposed the term ``differential equation$"$ (\ae quatio differentialis) and Leonhard Euler  in Saint-Petersburg obtained the first existence theorems \cite{Euler} using (in modern terms) the method of approximating the solution by a convergent sequence of piece-wise linear functions, the method that bears Euler's name today.

 Vyacheslav Alexandrovich Dobrovolskii \cite{Dobrovolskiy} describes the history of the Existence and Uniqueness Theorem for ordinary differential equations in the XIX century. Dobrovolskii conjectures that Augustin-Louis Cauchy may have been aware of the Existence and Uniqueness Theorem already by 1820\footnote{Cauchy's  notes of his course at the \'Ecole royale polytechnique, published by Charles Gilain \cite{Gilain}, are discussed by Adolf Pavlovich Yushkevitch in \cite{Yushk}}. In 1844 the Abbott Fran\c{c}ois-Marie-Napol{\'e}on Moigno published a course on in integral calculus based on Cauchy's lectures and included a proof of existence, using Euler's method, in lessons 26--27 (page 385 et seq). Moigno gives the sufficient condition that the derivative $\frac{\partial f(x, y)}{\partial y}$ be continuous and bounded for the existence of solutions of the ordinary differential equation written in the form
\[dy = f(x, y) dx\]
 Uniqueness of solutions is very briefly  mentioned. The extension to a system of ordinary differential equation is treated by Moigno in Lesson 39 (page 513 et seq). Again the uniqueness of solutions is treated very briefly.
 
The Existence and Uniqueness Theorem under the Lipschitz condition is due to Rudolf Lipschitz \cite{Lipschitz} who again uses the Euler method and cites neither Cauchy nor Moigno; Lipschitz gives a detailed proof for both the existence and uniqueness of solutions and makes the remark that uniqueness need not hold if the Lipschitz condition is not verified --- without, however, providing specific counterexamples.

We now proceed to a brief review of iteration proofs of the Existence and Uniqueness Theorem.
Let us first recall that solving equations by iterating a contraction is a procedure with counterparts in Antiquity. Heron of Alexandria in his \textit{Metrica} (AD 69) uses, in modern terms, the contraction
\[x \to \frac{1}{2}\left(x + \frac{R}{x}\right)\tag{13}\label{modern term}\]
in order to find the approximate value of $\sqrt{R}$. The very precise approximation to $\sqrt2$ on the tablet YBC 7289, palaeographically  dated to ``the first third of second millennium BC$"$ \cite{YBC 7289} and currently owned by the Yale Babylonian Collection, has led scholars to the conjecture that the scheme (\ref{modern term}) was already in use during the Old Babylonian Period.

 For solutions of ordinary differential equations an iterative scheme foreshadowing the Picard method already appears, as Dobrovolskii indicates, in the paper \cite{Dobrovolskiy} by Joseph Liouville, one of the first contributions to what later came to the called the Sturm-Lioville theory. \'{E}mile Picard first introduced the iterative method that now bears his name in the study of the partial differential equation:
\[\frac{\partial^2 u}{\partial x^2} + \frac{\partial^2 u}{\partial y^2} + a\frac{\partial u}{\partial x} + b\frac{\partial u}{\partial x} + cu = 0.\]
\'{E}mile Picard applied the method to the study of ordinary differential equations in a subsequent paper \cite{Dobrovolskiy}. Ernst Lindel\"{o}f \cite{Lindelof} gives the proof essentially in modern form, together with the estimate $O(\alpha^n/n!)$ for the rate of convergence (cf.(5), p. 121 in \cite{Lindelof}). Note, however, that Lindel\"{o}f includes proof of uniqeness of solutions by contradiction with the Lipschitz condition, and does not seem to include in his analysis the observation that the Picard iteration scheme does, by itself, furnish both the proof of existence and that of uniqueness.

We now proceed to a brief review of topological proofs of the Existence and Uniqueness Theorem.
A particular case of what came to be called the Brouwer Fixed Point Theorem was obtained by Piers Georgievich Bohl in 1904\cite{Bohl}. The results of Bohl are specifically geared towards applications to ordinary differential equations --- ``k\"{o}nnen dazu dienen, perodische L\"{o}sungen von Differentialgleichungen nachzuweisen$"$\footnote{can be used to obtain perodic solutions of differential equations} (p.186 in \cite{Bohl}; see a brief discussion by Myshkis and Rabinovich in \cite{Myshkis Rabinovich}).

 Let us now turn to the XX century developments in Lvov following the account of R.Duda \cite{Duda}. Stefan Banach established the Contraction Mapping Principle in his 1920 Ph.D. thesis published in 1922 \cite{Banach}. The Contraction Mapping Principle is formulated by Banach for self-maps of complete normed vector spaces --- the spaces that are called  Banach spaces today and that are, precisely, introduced in Banach's Ph.D. thesis.

 In \cite{Banach} Stefan Banach did not discuss possible applications of his Principle to ordinary differential equations. A few years later, however, the idea of proving the exstence of solutions for differential equations using topological considerations based on the Brouwer Fixed Point Theorem \cite{Brouwer} started to flourish in the Lvov School. Juliusz Schauder \cite{Schauder}, building on previous work of Birkhoff and Kellogg \cite{Birkhoff-Kellogg}, proved his famous Fixed Point Theorem precisely with a view towards solving differential equations. Let us emphasize that the topological methods of Schauder only give existence of solutions, not their uniqueness, much less a procedure for finding approximate solutions. In the USSR, the investigations of Schauder were pursued by A.N. Tikhonov \cite{Tikhonov}, who, again, only considered the existence of solutions. The formulation of the Contraction Mapping Principle for general metric spaces, and its application to the proof of existence and uniqueness of solutions of ordinary differential equations, is due to Renato Caccioppoli.

 Renato Caccioppoli is the tragic hero of XX century Italian mathematics.
Caccioppoli made a determinating contribution to the development of the Italian school of analysis: in the words of Alessandro Fig\`{a} Talamanca\footnote{Alessandro Fig{\`a} Talamanca, Renato Caccioppoli, Enciclopedia italiana, 2004.}, Caccioppoli's ideas ``ebbero un'influenza decisiva sullo sviluppo della analisi matematica in Italia: si può dire che un'intera generazione di analisti ha attinto da lui nuove idee, indirizzi e ispirazioni in un periodo in cui l'Italia si era isolata culturalmente dal resto del mondo (…)si deve pertanto a lui se, nel dopoguerra, l'analisi matematica in Italia ha potuto reinserirsi (…) nelle grandi correnti del pensiero matematico mondiale$"$\footnote{had a decisive influence on the development of mathematical analysis in Italy: it can be said that an entire generation of analysts obtained from him their new ideas, direction and inspiration in a period in which Italy had culturally isolated itself from the rest of the world (...)it is therefore due to him that, in the post-war period, mathematical analysis in Italy was able to reinsert itself ( ... ) in the great currents of world mathematical thought}.

 The years between the unification of Italy and the First World War are often defined as ``the golden age$"$ of Italian mathematics. The buoyant spirit of this optimistic period can be seen from the famous words of Eugenio Beltrami to Ernesto Ces\`{a}ro\footnote{Luigi Bianchi, Eugenio Beltrami, Enciclopedia italiana, 1930}, in which the Lombard mathematician, himself son of a miniaturist and grandson of a gem engraver, compliments the book of his Parthenopean colleague for ``davvero il requisito dell’italianit{\`a}, vale a dire di quel {\it quid} che risulta dal connubio della seriet{\`a} coll’agilit{\`a} della parola e del pensiero, cio{\`e} dell’elaborazione artistica del materiale scientifico$"$\footnote{really the requirement of Italianness, that is to say of that quid that results from the union of seriousness with the agility of speech and thought, that is, of the artistic elaboration of scientific material} --- words that Vito Volterra \cite{Volterra} describes as a sober, efficient and precise characterization of ``la produzione matematica italiana non solo recente ma di tutti i tempi$"$\footnote{Italian mathematical production not only recent but of all time}.

An avid reader with a deep interest for Proust and a special affinity for Rimbaud, whose portrait, together with that of Galois, he kept on his working table, an organizer, in Naples, of a cinema club devoted to contemporary film, a talented pianist who had contemplated a professional performing career, would Caccioppoli have recognized himself in Beltrami's words? Ennio de Giorgi says: ``per quanto sia difficile e incauto entrare nel mistero di un uomo, se dovessi vedere un filo tra l’interesse artistico, l’interesse scientifico, l’interesse sociale e civile di Caccioppoli, lo vedrei in questa aspirazione di fondo all’armonia, e nel dolore che tutte le varie disarmonie ai vari livelli gli procuravano” and develops “Perché c’è in lui (…) l’idea (…) dell’armonia pitagorica(…) che (…) la costruzione matematica (…) deve essere (…) bella e armonica, non può essere  (…) priva di bellezza$"$\footnote{as difficult and imprudent as it is to try to enter into the mystery of a man, if I were to look for a thread joining the artistic interest, the scientific interest, the social and civil interest of Caccioppoli, I would see it in his profound yearning for harmony, and in the pain that all the various disharmonies at various levels gave him (...) Because there is in him ( ... ) the idea ( ... ) of Pythagorean harmony ( ... ) that ( ... ) a mathematical construction ( ... ) must be ( ... ) beautiful and harmonious, it cannot be ( ... ) devoid of beauty}.\\
Entering the mathematical stage as the golden age of Italian mathematics had just vanished in the flames of the First World War, Renato Caccioppoli would himself live in a very different world from that of Eugenio Beltrami.

Renato Caccioppoli was born in Naples on 20 January 1904 into the family of the surgeon Giuseppe Caccioppoli, private physician to the Queen Margherita, and Sofia Bakunina, the daughter of Mikhail Alexandrovich Bakunin\footnote{Renato Caccioppoli left a profound impression. Several surveys, in particular, those published on his centenary in 2004,  are devoted to Caccioppoli's Mathematics. Five biographies as well as the  film  <<Morte di un matematico napoletano>>(1992), directed by Mario Martone and loosely based on the events of the mathematician's life, present Caccioppoli's person to the general public. In these brief remarks we follow the survey by Carlo Sbordone \cite{sbordone} and the biography by Lorenza Foschini \cite{foschini-it}}. On obtaining maturità classica, Caccioppoli entered the University of Naples to study engineering, but then switched to Mathematics. In the last years of the XIX century Nicola Trudi, Gabriele Torelli (whose son Ruggiero proved what is today called the Torelli Theorem), Ernesto Ces{\`ar}o, Alfredo Capelli all worked at the University of Naples. During Caccioppoli's years as a student, Pasquale del Pezzo, Duke of Cajanello, held the Chair of Projective Geometry at the University, whose Rector he had been from 1919 to 1921\cite{DelPezzo}. Of decisive impact on Caccioppoli was the arrival to Naples of Mauro Picone, whose assistant Caccioppoli became. Caccioppoli's active public stance made his life in fascist Italy difficult and dangerous\footnote{Caccioppoli's meeting with Andr\'e Gide in Sorrento in 1937 — a meeting that Gide in a December 1945 journal entry described as ``unforgettable$"$ — gave rise to a report by the secret police: ``Divisione polizia politica. Roma 27 agosto 1937.(...) Il Caccioppoli attribuisce a Gide queste parole: $"$Un regime, che poggia soltanto sulla fortuna personale di un uomo, non pu\`o essere che un regime transitorio! [Political police division. Rome 27 August 1937.(...) Caccioppoli attributes these words to Gide: $"$A regime, which rests only on the personal fortune of a man, can only be a transitional regime!$"$]}. In 1938 in order to avert a disaster the family had to place him into a psychiatric asylum; his piano was brought to him so that he could continue to practise. The fascist regime collapsed, but in the post-war Democratic-Christian Italy Caccioppoli remained suspect, unable, in particular, to travel abroad: with sneering brutality, his exit visas, while granted on paper, were issued in such a manner as to render actual travel impossible. As just one example, Caccioppoli was even prevented from going to Amsterdam in 1954. Facing the abyss of humiliation proved unbearable. In Caccioppoli's own words: $"$Napoli \`e una palude e noi siamo la fauna malata di questa palude. La vigliaccheria ci fa ingrassare e ci uccide contemporaneamente$"$ (Naples is a swamp and we are the sick fauna of this swamp. Cravenness makes us fat and kills us at the same time). Renato Caccioppoli took his life on 8 May 1959.

 Caccioppoli's work on elliptic equations developed, to a degree, in rivality with that of Petrovsky, but we are not aware of any direct contacts between Caccioppoli and the mathematicians of the Moscow School (cf. Caccioppoli's non-meeting with Kolmogorov in 1954). Nonetheless, a few years after the appearance, during his feverishly productive period, of the 1930 note on existence and uniqueness of solutions to differential equations Caccioppoli's work was included in the 1939 textbook, based on the author's 1936 Saratov lectures, by Ivan Georgievich Petrovsky.

\begin{center}
	\includegraphics[scale=0.2]{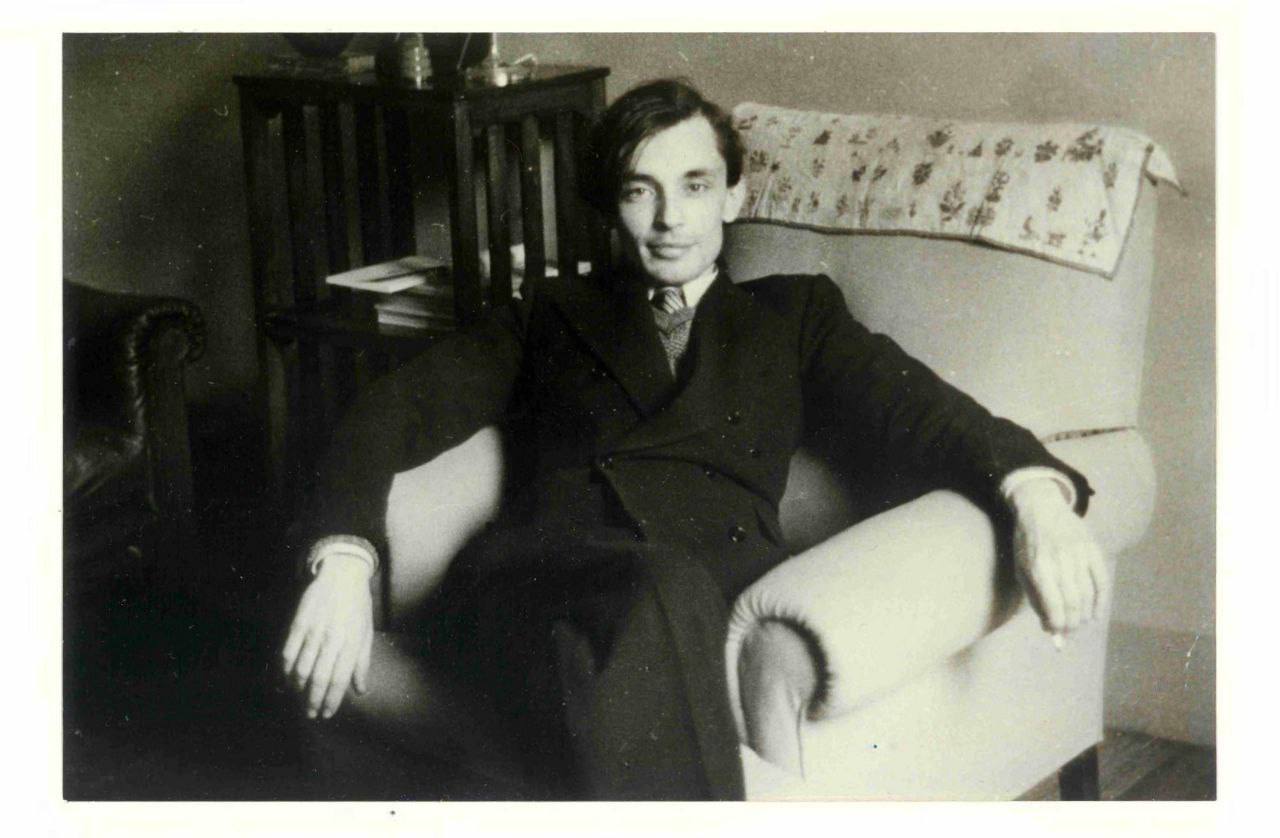}\\
	Renato Caccioppoli at the Levi-Civita family home. Photograph from the Levi-Civita archive, Ceccherini-Silberstein family.
\end{center}

\section*{ACKNOWLEDGEMENTS}
We are deeply grateful to the heirs of Tullio Levi-Civita for the kind permission to use the photograph of Renato Caccioppoli  in this paper. We are deeply grateful to the anonymous referees for a careful reeding of the manuscript and thoughtful and constructive criticism. We are deeply grateful to Yulij Sergeevich Ilyashenko for many helpful suggestions and, in particular, for indicating Dobrovolskii’s book to us. We are deeply grateful to Alexey Vladimirovich Klimenko for pointing out the evolution of nomenclature in different editions of Petrovsky’s classical textbook, to Tullio Ceccherini-Silberstein and Mauro Mariani for many helpful discussions, and to Michail Antonovich Borovikov for helpful remarks.
The authors are winners of the “BASIS” Foundation Competition in Mathematics and theoretical physics and are deeply grateful to the Jury
and the sponsors.

\section*{FUNDING}
This work was supported by the Russian Science Foundation under grant no. 24-71-10109, \\https://rscf.ru/en/project/24-71-10109/.






\phantom{123}\\
Alexander I. Bufetov\\
Steklov Mathematical Institute of Russian Academy of Sciences\\
8 Gubkina St., Moscow 119991, Russia\\
E-mail: \texttt{bufetov@mi-ras.ru}\\
Ilya I. Zavolokin\\
Steklov Mathematical Institute of Russian Academy of Sciences,\\
8 Gubkina St., Moscow 119991, Russia\\
E-mail: \texttt{tunderavatar@yandex.ru}\\
\end{document}